\documentclass[12pt,english]{amsart}
\usepackage[T1]{fontenc}
\usepackage[latin9]{inputenc}
\usepackage{fancyhdr}
\usepackage{amsthm}
\usepackage{amssymb}
\usepackage{esint}

\numberwithin{equation}{section} %% Comment out for sequentially-numbered
\numberwithin{figure}{section} %% Comment out for sequentially-numbered
\theoremstyle{plain}
\theoremstyle{plain}
\newtheorem{thm}{Theorem}
  \theoremstyle{remark}
  \newtheorem*{acknowledgement*}{Acknowledgement}
  \theoremstyle{plain}
  \newtheorem{lem}[thm]{Lemma}
  \theoremstyle{plain}
  \newtheorem{cor}[thm]{Corollary}
  \theoremstyle{remark}
  \newtheorem*{rem*}{Remark}

\usepackage{babel}

\begin{document}
\global\long\def\A{{\mathcal{A}}}

\global\long\def\B{{\mathcal{B}}}

\global\long\def\C{{\mathcal{C}}}

\global\long\def\M{{\mathcal{M}}}

\global\long\def\N{{\mathcal{N}}}

\global\long\def\O{{\mathcal{O}}}

\global\long\def\bbc{{\mathbb{C}}}

\global\long\def\bbp{{\mathbb{P}}}

\global\long\def\bbd{{\mathbb{D}}}

\global\long\def\bbt{{\mathbb{T}}}

\global\long\def\Bbr{{\mathbb{R}}}

\global\long\def\Bbn{{\mathbb{N}}}

\global\long\def\dprime{{\prime\prime}}

\global\long\def\vol#1{\mathrm{vol}(#1)}

\global\long\def\volcn#1{\mathrm{vol_{\bbc^{N}}}(#1)}

\global\long\def\volrn#1#2{\mathrm{vol_{\Bbr^{#1}}}(#2)}

% complex analysis

\global\long\def\disk#1{#1 \bbd}

\global\long\def\cir#1{#1 \bbt}

% probability

\global\long\def\pr#1{\bbp\left( #1 \right)}

\global\long\def\Ex{\mathbb{E}}

\global\long\def\Cov#1{\mathrm{Cov}(#1)}

% special functions

\global\long\def\Li{{\mathrm{Li}}}

\global\long\def\erf{\mathrm{erf}}

\global\long\def\erfc{\mathrm{erfc}}

\title{Asymptotics of The Hole Probability for Zeros of Random Entire Functions}

\author{Alon Nishry}
\begin{abstract}
Consider the random entire function\[
f(z)=\sum_{n=0}^{\infty}\phi_{n}\,\frac{z^{n}}{\sqrt{n!}},\eqno(*)\]
 where the $\phi_{n}$ are i.i.d. standard complex Gaussian variables.
The zero set of this function is distinguished by invariance of its
distribution with respect to the isometries of the plane.

We study the probability $P_{H}(r)$ that $f$ has no zeroes in the
disk $\left\{ |z|<r\right\} $ (hole probability). Improving a result
of Sodin and Tsirelson, we show that\[
\log P_{H}(r)=-\frac{e^{2}}{4}\cdot r^{4}+o(r^{4})\]
as $r\to\infty.$ The proof does not use distribution invariance of
the zeros, and can be extended to other Gaussian Taylor series.

If instead of Gaussians we take Rademacher or Steinhaus random variables
$\phi_{n}$, we get a very different result. There exists $r_{0}$
so that every random function of the form ($*$) with Rademacher or
Steinhaus coefficients must vanish in the disk $\{|z|<r_{0}\}$. 
\end{abstract}
\maketitle

\section{Introduction}

Consider the following random entire function

\begin{equation}
f(z)=\sum_{n=0}^{\infty}\phi_{n}a_{n}z^{n},\label{eq:f_definition}\end{equation}
where $a_{n}=\left(n!\right)^{-1/2}$ and $\phi_{n}$ are independent
standard complex Gaussian random variables (i.e., each $\phi_{n}$
has the density function $\pi^{-1}\cdot\exp\left(-|z|^{2}\right)$
with respect to Lebesgue measure on $\bbc)$. The random zero set
of this function is known to be distribution invariant with respect
to isometries of the plane. Furthermore, this is the only Gaussian
entire function with distribution invariant zeros (see \cite{ST1},
and the forthcoming book \cite{BKPV} for details and discussion).

One of the interesting characteristics of the random zero process
$f^{-1}\left\{ 0\right\} $ is the asymptotic decay of the event where
$f$ has no zeros inside the disk $\left\{ |z|\le r\right\} $ when
$r$ is large. Since the decay rate is known to be exponential, we
use the notation\[
p_{H}(z)=\log^{-}P_{H}(r)=\log^{-}\pr{f(z)\ne0\mbox{ for }|z|\le r}.\]
In the paper \cite{ST3}, Sodin and Tsirelson showed that for $r\ge1$,\begin{equation}
c_{1}r^{4}\le p_{H}(r)\le c_{2}r^{4}\label{eq:hole_ST_asymp}\end{equation}
with some positive numerical constants $c_{1}$ and $c_{2}$. This
result was extended in different directions by Ben Hough \cite{BH},
Krishnapur \cite{K}, Zrebiec \cite{Zr1,Zr2} and Shiffman, Zelditch
and Zrebiec \cite{SZZ}.

In \cite{ST3}, Sodin and Tsirelson raised the question whether the
limit\[
\lim_{r\to\infty}\frac{p_{H}(r)}{r^{4}}\]
exists and what is its value ? Our main result answers this question
(and estimates the remainder):
\begin{thm}
For $r$ large enough\label{thm:Main_Thm_Intro}\begin{equation}
p_{H}(r)=\frac{e^{2}}{4}\cdot r^{4}+\O\left(r^{18/5}\right).\label{eq:accurate_hole_asymp}\end{equation}

\end{thm}
The constant $\frac{e^{2}}{4}$ arrives as follows. We introduce
the function\[
S(r)=\log\prod_{\left\{ n\,:\, a_{n}r^{n}\ge1\right\} }\left(a_{n}r^{n}\right)^{2}=2\cdot\sum_{\left\{ n\,:\, a_{n}r^{n}\ge1\right\} }\log\left(a_{n}r^{n}\right),\]
and prove that\begin{equation}
p_{H}(r)=S(r)+\O\left(r^{18/5}\right),\qquad r\to\infty.\label{eq:hole_prob_wrt_S_r}\end{equation}
Then it is easy to see that\[
S(r)=\frac{e^{2}}{4}\cdot r^{4}+\O\left(r^{2}\log r\right),\qquad r\to\infty.\]
Actually, it is plausible that estimate \eqref{eq:accurate_hole_asymp}
holds with a better estimate of the remainder, for instance, with
$\O\left(r^{2+\epsilon}\right)$ with any $\epsilon>0$.

The proof of the upper bound in \eqref{eq:hole_prob_wrt_S_r} is similar
to the proof of the upper bound in \eqref{eq:hole_ST_asymp} given
in \cite{ST3}, the only difference is that our estimates are slightly
more accurate (note that in \cite{ST3} this is considered to be the
lower bound). The proof of the lower bound in \eqref{eq:hole_prob_wrt_S_r}
combines techniques from \cite{ST3} with direct estimates of probability
of some events in high-dimensional linear spaces with Gaussian measure.
Note that somewhat similar ideas were used in \cite{SZZ}.

In contrast to \cite{ST3}, our proof of the lower bound in \eqref{eq:hole_prob_wrt_S_r}
does not use distribution invariance of the zero set of $f$, and
our main result can be extended to other Gaussian entire functions
of the form \eqref{eq:f_definition} with regular sequence of the
coefficients $a_{n}$. For instance, one may consider Gaussian Mittag-Leffler
functions\[
f(z)=\sum_{n=0}^{\infty}\phi_{n}\cdot\frac{z^{n}}{\Gamma\left(\alpha n+1\right)}\]
with $\alpha>0$. In this case, the corresponding function $S(r)$
has the asymptotics\[
S(r)=\frac{1}{2\alpha}r^{2/\alpha}+\O_{\alpha}\left(r^{1/\alpha}\log r\right)\]
and then only minor modifications in the proof of Theorem \ref{thm:Main_Thm_Intro}
are needed to show that\[
p_{H}(r)=S(r)+\O_{\alpha}\left(r^{9/5\alpha}\right).\]

One can ask what happens when the i.i.d. coefficients $\zeta_{n}$
in \eqref{eq:f_definition} are not Gaussian ? The following deterministic
result shows that the situation might be very different.
\begin{thm}
Let $K\subset\bbc$ be a compact set and $0\notin K$. Suppose that
$\mathbf{\phi}_{n}\in K$ for each $n$, and that $a_{n}=\left(n!\right)^{-1/2}.$
Then there exists $r_{0}(K)<\infty$ so that $f(z)$ must vanish somewhere
in the disk $\left\{ |z|\le r_{0}(K)\right\} $.\label{thm:Compact_Coeffs_Thm}
\end{thm}
The idea of the proof is very simple. A standard compactness argument
shows that if the result is wrong then there exists an entire function
$f$ of the form \eqref{eq:f_definition} with $\phi_{n}\in K$ and
$a_{n}=\left(n!\right)^{-1/2}$ that does not vanish on $\bbc$. Since
$f$ is an entire function of order $2$, it equals $\exp\left(\alpha z^{2}+\beta z+\gamma\right)$
with complex constants $\alpha,\beta,\gamma$. Then it is not difficult
to verify that the Taylor coefficients of that function cannot be
equal to $\phi_{n}/\sqrt{n!}$ with $\zeta_{n}\in K$.
\begin{acknowledgement*}
First of all I would like to thank my advisor Mikhail Sodin for his
more than substantial involvement in writing this paper, and supporting
me in general. Starting with the subject, continuing with technical
issues, and especially in the presentation, his fingerprints are found
everywhere in this paper. I would also like to thank Fedor Nazarov
for contributing his idea regarding the second theorem.
\end{acknowledgement*}

\section{Notations and Elementary Estimates\label{sec:Elementary-Estimates}}

In what follows we frequently use that if $w$ is a standard Gaussian
random variable, then, \begin{equation}
\pr{|w|\ge\lambda}=\exp(-\lambda^{2}),\label{eq:prob_large}\end{equation}
 and for \[
\lambda\le1,\]

\begin{equation}
\pr{|w|\le\lambda}\in\left[\frac{\lambda^{2}}{2},\lambda^{2}\right].\label{eq:prob_small}\end{equation}
We denote by $\disk r$the disk $\left\{ z\,:\,|z|<r\right\} $ and
by $\cir r$ its boundary $\left\{ z\,:\,|z|=r\right\} $. The letter
$C$ denotes positive numerical constant (which can change between
lines).

In what follows, we use several elementary estimates, we skip their
proofs.
\begin{lem}
The sequence $a_{n}r^{n}$ has a local maximum only in the interval
$n\in\left\{ \left\lceil r^{2}-1\right\rceil ,\left\lfloor r^{2}\right\rfloor \right\} $.\label{lem:a_n_r^n_unimodal}
\end{lem}
Using Striling's approximation we have
\begin{lem}
\label{lem:Basic-Estimates-1}For all $n\ge1$ we have\begin{equation}
\frac{1}{\sqrt{3n}}\left(\frac{e}{n}\right)^{\frac{n}{2}}\le a_{n}\le\left(\frac{e}{n}\right)^{\frac{n}{2}},\label{eq:c_n_bounds}\end{equation}
moreover for $n\in\left\{ 1,\ldots,\left\lfloor er^{2}\right\rfloor \right\} $
\begin{equation}
\frac{1}{3r}\cdot(a_{n}r^{n})^{-1}\le1\label{eq:mid_range_bound}\end{equation}
and for $n\ge er^{2}$\begin{equation}
\log\left(a_{n}r^{n}\right)\le-\frac{1}{2}\left(n-er^{2}\right).\label{eq:h_n_tail}\end{equation}

\end{lem}
Since by Lemma \ref{lem:a_n_r^n_unimodal} we see the sequence $a_{n}r^{n}$
is unimodal in $\left[0,\left\lfloor er^{2}\right\rfloor \right]$,
it is easy to estimate $S(r)$ with a corresponding integral and get
\begin{lem}
We have\[
S(r)=\frac{e^{2}}{4}\cdot r^{4}+\O\left(r^{2}\log r\right).\]

\end{lem}

\section{Upper Bound for $p_{H}(r)$\label{sec:Thm_1-Upper-Bound}}

In this section, we show that for $r$ large enough,

\[
p_{H}(r)\le S(r)+\O\left(r^{2}\log r\right).\]

\begin{proof}
Denote by $\Omega_{r}$ the following event\[
\begin{array}{lll}
\mbox{{\rm \mbox{(i)}}} & |\phi_{0}|\ge2r & ,\ \\
\mbox{{\rm \mbox{(ii)}}} & |\phi_{n}|\le\frac{1}{3r}\cdot(a_{n}r^{n})^{-1} & n\in\left\{ 1,\ldots,\left\lfloor er^{2}\right\rfloor \right\} ,\\
\mbox{{\rm \mbox{(iii)}}} & |\phi_{n}|\le\exp\left(\frac{n-er^{2}}{4}\right) & n\ge\left\lfloor er^{2}\right\rfloor +1.\end{array}\]
We prove that if $r$ is large enough and the event $\Omega_{r}$
occurs, then $f(z)\ne0$ inside $\disk r$, and that\[
\log\pr{\Omega_{r}}\ge-S(r)-\O\left(r^{2}\log r\right).\]
Note that\begin{equation}
|f(z)|\ge|\phi_{0}|-\sum_{n=1}^{\infty}|\phi_{n}|a_{n}r^{n}.\label{eq:Gen_1}\end{equation}
First, estimate the sum \[
\sum_{n=1}^{\left\lfloor er^{2}\right\rfloor }|\phi_{n}|a_{n}r^{n}\le\sum_{n=1}^{\left\lfloor er^{2}\right\rfloor }\frac{1}{3r}\le r.\]
To bound the tail we use \eqref{eq:h_n_tail}, \begin{eqnarray*}
\sum_{n\ge\left\lfloor er^{2}\right\rfloor +1}|\phi_{n}|a_{n}r^{n} & \le & \sum_{\left\lfloor er^{2}\right\rfloor +1}^{\infty}\exp\left(\frac{n-er^{2}}{4}-\frac{1}{2}\left(n-er^{2}\right)\right)\\
 & \le & \sum_{n=0}^{\infty}\exp\left(-n/2\right)=\O(1).\end{eqnarray*}
From \eqref{eq:Gen_1}, we have\[
|f(z)|\ge2r-r-\O(1)>0,\]
 for $r$ large enough. We see that $f(z)\ne0$ inside $\disk r$.

Now we estimate the probability of $\Omega_{r}$ using \eqref{eq:prob_large}
and \eqref{eq:prob_small}. First,\[
\pr{\mbox{(i)}}=\exp\left(-4r^{2}\right).\]
For $n\ge\left\lfloor er^{2}\right\rfloor +1$, we have\[
\pr{\mbox{(iii)}_{n}}=1-\exp\left(-\exp\left(\frac{n-er^{2}}{2}\right)\right).\]
That is,\[
\pr{\mbox{(iii)}}=\prod_{n\ge\left\lfloor er^{2}\right\rfloor +1}\pr{\mbox{(iii)}_{n}}=\exp\left(\sum_{n\ge\left\lfloor er^{2}\right\rfloor +1}\log\pr{\mbox{(iii)}_{n}}\right).\]
Taking logarithm of $\pr{\mbox{(iii)}_{n}}$, we can see that we have
the following estimate for $r$ large enough,\[
\log\pr{\mbox{(iii)}_{n}}\ge-\exp\left(-\frac{n-er^{2}}{2}\right),\]
so $\pr{\mbox{(iii)}}$ is larger than some constant, which does not
depend on $r$. For the term $\pr{\mbox{(ii)}}$, recalling \eqref{eq:mid_range_bound},
we use the estimate\[
\pr{\mbox{(ii)}_{n}}\ge\frac{(a_{n}r^{n})^{-2}}{18r^{2}},\]
and get\[
\pr{\mbox{(ii)}}\ge\prod_{n=0}^{\left\lfloor er^{2}\right\rfloor }\frac{(a_{n}r^{n})^{-2}}{18r^{2}}=\exp\left(-S(r)-\left\lfloor er^{2}\right\rfloor \cdot\log18r^{2}\right).\]
Since $\pr{\Omega_{r}}=\pr{\mbox{(i)}}\pr{\mbox{(ii)}}\pr{\mbox{(iii)}}$,
we get the required result:\[
p_{H}(r)\le-\log\pr{\Omega_{r}}\le S(r)+\O\left(r^{2}\log r\right).\]

\end{proof}

\section{Lower Bound for $p_{H}(r)$\label{sec:Thm-1-Lower-Bound}}

In this section we show that for $r$ large enough\[
p_{H}(r)\ge S(r)-Cr^{18/5}.\]
Define $M(r)={\displaystyle \max_{|z|\le r}|f(z)|}$, we start by
studying the deviations of $\log M(r)$ from the mean $\frac{1}{2}r^{2}$.
Then we consider large deviations of the expression ${\displaystyle \intop_{\cir r}\log|f(z)|\, dm}$,
where $m$ is the normalized angular measure on $\cir r$. Finally,
we use the fact that if $n(r)=0$ then $\log|f(z)|$ is a harmonic
function inside $\disk r$ to get the result.

\subsection{Large deviations for $\log M(r)$}

We use the first part of Lemma 1 in the paper \cite{ST3} as
\begin{lem}
\label{lem:log_M_upp_bnd}Given $\sigma>0$, we have for $r$ large
enough\[
\log\pr{\frac{\log M(r)}{\frac{1}{2}r^{2}}\ge1+\sigma}\le-\exp\left(\sigma r^{2}\right).\]

\end{lem}
In the other direction we have
\begin{lem}
\label{lem:log_M_low_bnd}We have the following estimate for the lower
bound of $M(r)$\[
\log\pr{M(r)\le1}\le-S(r).\]
\end{lem}
\begin{proof}
Suppose $\log|f(z)|\le0$ in $\disk r$, then using Cauchy's estimate
for the coefficients of $f(z)$ we can get an estimate to the probability
of this event. We have\[
|\phi_{n}|a_{n}r^{n}\le M(r)\le1\]
or\[
|\phi_{n}|\le\left(a_{n}r^{n}\right)^{-1}.\]
The probability of each event, for $n\in\left\{ 0,\ldots,\left\lfloor er^{2}\right\rfloor \right\} $
is bounded by (using Lemma \ref{lem:Basic-Estimates-1})\[
\pr{|\phi_{n}|\le\left(a_{n}r^{n}\right)^{-1}}\le\left(a_{n}r^{n}\right)^{-2}.\]
We get\[
\pr{M(r)\le1}\le\prod_{n=0}^{\left\lfloor er^{2}\right\rfloor }\left(a_{n}r^{n}\right)^{-2}=\exp\left(-S(r)\right).\]

\end{proof}

\subsection{Discretization of the logarithmic integral}

In this section $\delta\in(0,1)$, $N=\left\lfloor er^{2}\right\rfloor $
, $\kappa=1-\delta^{1/2}$ and the points $\left\{ z_{j}\right\} _{j=0}^{N-1}$
are equally distributed on $\cir{\kappa r}$, that is \[
z_{j}=\kappa r\exp\left(\frac{2\pi ij}{N}\right).\]
Also $m$ is the normalized angular measure on $\cir r$. Under this
conditions we have
\begin{lem}
\label{lem:approx_log_int}Outside an exceptional set of probability
at most\[
2\exp\left(-S(\kappa r)\right)\]
we have\begin{equation}
\frac{1}{N}\sum_{j=0}^{N-1}\log|f(z_{j})|\le\intop_{\cir r}\log|f|\, dm+\frac{C}{\delta^{2}}\label{eq:error_in_discrete_approx}\end{equation}
with $C$ a positive numerical constant.\end{lem}
\begin{proof}
Denote by $P_{j}(z)=P(z,z_{j})$ the Poisson kernel for the disk $\disk r,$$|z|=r$,
$|z_{j}|<r$. Since $\log|f|$ is a subharmonic function we have\begin{eqnarray*}
\frac{1}{N}\sum_{j=0}^{N-1}\log|f(z_{j})| & \le & \intop_{\cir r}\left(\frac{1}{N}\sum_{j=0}^{N-1}P_{j}\right)\log|f|\, dm\\
 & = & \intop_{\cir r}\log|f|\, dm+\intop_{\cir r}\left(\frac{1}{N}\sum_{j=0}^{N-1}P_{j}-1\right)\log|f|\, dm.\end{eqnarray*}
The last expression can be estimated by \begin{equation}
\intop_{\cir r}\left(\frac{1}{N}\sum_{j=0}^{N-1}P_{j}-1\right)\log|f|\, dm\le\max_{z\in\cir r}\left|\frac{1}{N}\sum_{j=0}^{N-1}P_{j}-1\right|\cdot\intop_{\cir r}\left|\log|f|\right|\, dm.\label{eq:log_int_error_terms}\end{equation}
For the first factor in the RHS of \eqref{eq:log_int_error_terms},
we start with\[
\intop_{\cir{\kappa r}}P(z,\omega)\, dm(\omega)=1,\]
and then split the circle $\cir{\kappa r}$ into a union of $N$ disjoint
arcs $I_{j}$ of equal angular measure $\mu(I_{j})=\frac{1}{N}$ centered
at the $z_{j}$'s. Then\[
1=\frac{1}{N}\sum_{j=0}^{N-1}P(z,z_{j})+\sum_{j=0}^{N-1}\intop_{I_{j}}\left(P(z,\omega)-P(z,z_{j})\right)\, dm(\omega),\]
and\begin{eqnarray*}
|P(z,\omega)-P(z,z_{j})| & \le & \max_{\omega\in I_{j}}|\omega-z_{j}|\cdot\max_{z,\omega}|\nabla_{\omega}P(z,\omega)|\\
 & \le & \frac{2\pi r}{N}\cdot\frac{Cr}{(r-|\omega|)^{2}}\le\frac{C}{\delta N}.\end{eqnarray*}
For the second factor on the RHS of \eqref{eq:log_int_error_terms},
using Lemma \ref{lem:log_M_low_bnd}, we may suppose that there is
a point $a\in\kappa\cir r$ such that $\log|f(a)|\ge0$ (discarding
an exceptional event of probability at most $\exp\left(-S(\kappa r)\right)$).
Then we have\[
0\le\intop_{\cir r}P(z,a)\log|f(z)|\, dm(z),\]
and hence\[
\intop_{\cir r}P(z,a)\log^{-}|f(z)|\, dm(z)\le\intop_{\cir r}P(z,a)\log^{+}|f(z)|\, dm(z).\]
For $|z|=r$ and $|a|=\kappa r$ we have,\[
\frac{\delta^{\frac{1}{2}}}{2}\le\frac{1-(1-\delta^{\frac{1}{2}})}{1+(1-\delta^{\frac{1}{2}})}\le P(z,a)\le\frac{1+(1-\delta^{\frac{1}{2}})}{1-(1-\delta^{\frac{1}{2}})}\le\frac{2}{\delta^{\frac{1}{2}}}.\]
By Lemma \ref{lem:log_M_upp_bnd}, outside a very small exception
set (of the order $\exp\left(-\exp\left(r^{2}\right)\right)$), we
have\[
\intop_{\cir r}\log^{+}|f|\, d\mu\le r^{2},\]
and therefore\[
\intop_{\cir r}\log^{-}|f|\, d\mu\le\frac{Cr^{2}}{\delta}.\]
Finally\[
\intop_{\cir r}\left|\log|f|\right|\, d\mu\le\frac{Cr^{2}}{\delta}.\]
Using the fact that $N=\left\lfloor er^{2}\right\rfloor $, we see
that the lemma is proved.
\end{proof}

\subsection{Deviation from the logarithmic integral}

If we use the notation $\zeta_{j}=f(z_{j})$, we know that the Gaussian
random variables $\left\{ \zeta_{j}\right\} _{j=0}^{N-1}$ have a
multivariate complex Gaussian distribution, with covariance matrix
$\Sigma$, where\begin{eqnarray*}
\Sigma_{ij} & = & \Cov{\zeta_{i},\zeta_{j}}=\Cov{f(z_{i}),f(z_{j})}\\
 & = & \Ex(f(z_{i})\overline{f(z_{j})})=\sum_{n=0}^{\infty}a_{n}^{2}z_{i}\bar{z_{j}}=\exp\left(z_{i}\bar{z_{j}}\right).\end{eqnarray*}
We also know that the density function of this distribution is\[
\zeta\mapsto\frac{1}{\pi^{n}\cdot\det\Sigma}\cdot\exp(-\zeta^{*}\Sigma^{-1}\zeta).\]
We introduce the sets\begin{equation}
\A^{\prime}=\left\{ \zeta\,:\,\prod_{j=0}^{N-1}|\zeta_{j}|\le\exp\left(4N\cdot\log r+C\delta^{-2}\cdot r^{2}\right)\right\} \label{eq:def_set_a_prime}\end{equation}
and\begin{equation}
\A=\left\{ \zeta\,:\,\zeta\in\A^{\prime}\mbox{ and }|\zeta_{j}|=|f(z_{j})|\le\exp\left(2r^{2}\right),\quad0\le j\le N-1\right\} \label{eq:def_set_a}\end{equation}
and denote by $\B$ the set where estimate \eqref{eq:error_in_discrete_approx}
in Lemma \ref{lem:approx_log_int} holds. Using these notations we
get the simple
\begin{lem}
\label{lem:log_int_deviations}\[
\pr{\int_{\cir r}\log\,|f(z)|\, d\mu\le4\log r}\le\pr{\A}+\pr{\B^{c}}+\pr{\A^{\prime}\backslash\A}.\]
\end{lem}
\begin{proof}
We start by discarding the exceptional set in Lemma \ref{lem:approx_log_int},
this adds the term $\pr{\B^{c}}$. Now we can assume that\[
\frac{1}{N}\sum_{j=0}^{N-1}\log|f(z_{j})|\le\intop_{\cir r}\log|f|\, dm+\frac{C}{\delta^{2}},\]
or\[
\prod_{j=0}^{N-1}|\zeta_{j}|\le\exp\left(N\cdot\intop_{\cir r}\log|f|\, dm+\frac{C}{\delta^{2}}\cdot N\right).\]
In terms of probabilities we can write\[
\pr{\int_{\cir r}\log\,|f(z)|\, dm\le4\log r}\le\pr{\B^{c}}+\pr{\A^{\prime}},\]
and since\[
\pr{\A^{\prime}}=\pr{\A}+\pr{\A^{\prime}\backslash\A},\]
we get the required result.
\end{proof}
Before we continue, we need two asymptotic estimates.
\begin{lem}
Let $\Sigma$ be the covariance matrix defined above. We have the
following estimate\[
\log\left(\det\Sigma\right)\ge S(\kappa r).\]
\end{lem}
\begin{proof}
Notice that we can represent $\Sigma$ in the following form\[
\Sigma=V\cdot V^{*}\]
where \[
V=\left(\begin{matrix}a_{0} & a_{1}\cdot z_{1} & \dots & a_{N}\cdot z_{1}^{N} & \dots\\
\vdots & \vdots & \vdots & \vdots & \dots\\
a_{0} & a_{1}\cdot z_{N} & \dots & a_{N}\cdot z_{N}^{N} & \dots\end{matrix}\right).\]
By the Cauchy-Binet formula we have \[
\det\ \Sigma=|\Sigma|=\sum_{t}|m_{t}(V)|^{2},\]
where the sum is taken over all principal minors $m_{t}(V)$ of the
matrix $V$.

In our case the following minor is sufficient for the estimates\begin{eqnarray*}
|\Sigma| & \ge & \left|\begin{array}{cccc}
a_{1}z_{1} & a_{2}z_{1}^{2} & \ldots & a_{N}z_{1}^{N}\\
\vdots & \vdots & \vdots & \vdots\\
a_{1}z_{N} & a_{2}z_{N}^{2} & \ldots & a_{N}z_{N}^{N}\end{array}\right|^{2}\\
 & = & \prod_{n=1}^{N}a_{n}^{2}\cdot\prod_{i=1}^{N}\left|z_{i}\right|^{2}\cdot\prod_{1\le i\ne j\le N}\left|z_{i}-z_{j}\right|\\
 & = & \Pi_{1}\cdot\Pi_{2}\cdot\Pi_{3}.\end{eqnarray*}
It is clear that\[
\Pi_{2}=\left(\kappa r\right)^{2N}.\]
The $z_{i}$'s are the roots of the equation $z^{N}=\left(\kappa r\right)^{N}$,
denoting $z_{1}=\kappa r$ we get\[
\prod_{i=2}^{N}(z_{1}-z_{i})=N\cdot\left(\kappa r\right)^{N-1},\]
and\[
\Pi_{3}=\prod_{1\le i\ne j\le N}|z_{i}-z_{j}|=\left(\prod_{i=2}^{N}\left|z_{1}-z_{i}\right|\right)^{N}=\left(\kappa r\right)^{N(N-1)}\cdot N^{N}.\]
We now {}``collect'' the product of the $\kappa r$'s and rewrite
it as $\left(\kappa r\right)^{N(N-1)}=\prod_{n=1}^{N-1}\left(\kappa r\right)^{2n}$
and get\[
\det\Sigma\ge\prod_{n=1}^{N}a_{n}^{2}\left(\kappa r\right)^{2n}.\]
Using the fact that $N=\left\lfloor er^{2}\right\rfloor $, we have\[
\det\Sigma\ge\exp\left(S(\kappa r)\right).\]

\end{proof}
We denote by $I$ the following quantity\begin{equation}
I=\pi^{-N}\cdot\volcn{\A}.\label{eq:Asymp-Int}\end{equation}
We use the following lemma to estimate $I$,
\begin{lem}
Set $s>0$, $t>0$ and $N\in\Bbn^{+}$, such that $\log\left(t^{N}/s\right)\ge N$.
Denote by $\C_{N}$ the following set\[
\C_{N}=\C_{N}\left(t,s\right)=\left\{ r=\left(r_{1},\ldots,r_{N}\right)\::\:0\le r_{j}\le t,\,\prod_{1}^{N}r_{j}\le s\right\} .\]
Then\[
\volrn N{\C_{N}}\le\frac{s}{\left(N-1\right)!}\log^{N}\left(t^{N}/s\right).\]
\end{lem}
\begin{proof}
We will find an expression for the volume using induction. First we
define for $k\ge1$ \[
V_{k}(t,s)=\volrn k{\C_{k}(t,s)}.\]
We notice that if $s\ge t^{k}$ then $V_{k}(t,s)=t^{k}.$ Now we can
write\[
V_{k}(t,s)=\intop_{0}^{t}\left(\intop_{0}^{t}\ldots\intop_{0}^{t}\chi\left(\prod_{2}^{k}r_{j}\le s/x\right)\, dr_{2}\ldots dr_{k}\right)\, dx,\]
where $\chi$ is the characteristic function of the set $\left\{ r\,:\,\prod_{2}^{k}r_{j}\le s/x\right\} $,
if $s<t^{k}$ we can rewrite this expression as\begin{eqnarray*}
V_{k}(t,s) & = & \intop_{0}^{a}V_{k-1}\left(t,\frac{s}{x}\right)\, dx\\
 & = & \intop_{0}^{s/t^{k-1}}V_{k-1}\left(t,\frac{s}{x}\right)\, dx+\intop_{s/t^{k-1}}^{t}V_{k-1}\left(t,\frac{s}{x}\right)\, dx=I_{1}+I_{2}.\end{eqnarray*}
For the first integral we have $x<s/t^{k-1}$ or $s/x>t^{k-1}$ and
so $V_{k-1}\left(t,\frac{s}{x}\right)=t^{k-1}$, meaning $I_{1}=s$.
We can now prove by induction\[
V_{k}(t,s)=\begin{cases}
t^{k} & \quad s\ge t^{k},\\
\sum_{m=0}^{k-1}\frac{s}{m!}\cdot\left(\log\left(t^{k}/s\right)\right)^{m} & \quad s<t^{k}.\end{cases}\]
For $k=1$ this is trivial, now using the expression above we have,
\begin{eqnarray*}
V_{k}(t,s) & = & s+\intop_{s/t^{k-1}}^{t}\sum_{m=0}^{k-2}\frac{s}{x\cdot m!}\cdot\left(\log\left(xt^{k-1}/s\right)\right)^{m}\, dx\\
 & = & s+\sum_{m=0}^{k-2}\left[\left.\frac{s}{\left(m+1\right)!}\log^{m+1}\left(\frac{t^{k-1}x}{s}\right)\right|_{x=s/t^{k-1}}^{t}\right]\\
 & = & \sum_{m=0}^{k-1}\frac{s}{m!}\cdot\left(\log\left(t^{k}/s\right)\right)^{m}.\end{eqnarray*}
We conclude that\[
\volrn N{\C_{N}}=\sum_{m=0}^{N-1}\frac{s}{m!}\cdot\left(\log\left(t^{N}/s\right)\right)^{m}.\]
Since $\log\left(t^{N}/s\right)\ge N$, we can approximate the integral
from above by\begin{eqnarray*}
V_{N}(t,s) & = & \sum_{m=0}^{N-1}\frac{s}{m!}\cdot\log^{m}\left(t^{N}/s\right)\le s\cdot\sum_{m=0}^{N-1}\frac{\log^{m}\left(t^{N}/s\right)}{m!}\cdot\frac{\log^{N-m}\left(t^{N}/s\right)}{(m+1)\cdot\ldots\cdot N}\\
 & \le & \frac{s}{\left(N-1\right)!}\log^{N}\left(t^{N}/s\right).\end{eqnarray*}
This finishes the proof.
\end{proof}
Now we have as an almost immediate
\begin{cor}
For $r$ large enough and $\delta\ge r^{2-\epsilon}$, we have\[
\log I\le C\left(\log r+\delta^{-2}\right)r^{2}.\]
\end{cor}
\begin{proof}
We recall that\[
\A=\left\{ \begin{array}{cc}
 & |\zeta_{j}|=|f(z_{j})|\le\exp\left(2r^{2}\right),\quad0\le j\le N-1\\
\zeta\,:\, & \mbox{and}\\
 & \prod_{j=0}^{N-1}|\zeta_{j}|\le\exp\left(4N\cdot\log r+C\delta^{-2}\cdot r^{2}\right)\end{array}\right\} .\]
To shorten the expressions we use $s=\exp\left(4N\cdot\log r+C\delta^{-2}\cdot r^{2}\right)$
and $t=\exp\left(2r^{2}\right)$. We notice that under the assumptions
that we made and for $r$ large enough $\log\left(t^{N}/s\right)\ge N$
(remember that $N=\left\lfloor er^{2}\right\rfloor $). We want to
translate the integral into an integral in $\Bbr^{N}$, using the
change of variables $\zeta_{j}=r_{j}\cos(\theta_{j})+ir_{j}\sin(\theta_{j})$.
Integrating out the variables $\theta_{j}$, we get $I^{\prime}=2^{N}\intop_{\C}\prod r_{j}\, dr$,
where the new domain is\[
\C=\left\{ r=\left(r_{1},\ldots,r_{N}\right)\::\:0\le r_{j}\le t,\,\prod_{j=1}^{N}r_{j}\le s\right\} .\]
We can find an explicit expression for this integral, but, instead
we will simplify it even more to\begin{equation}
I^{\prime}\le2^{N}\cdot s\cdot\volrn N{\C}\label{eq:I_prime_upp_bnd}\end{equation}
Now we can use the previous lemma, and get (for $r$ large enough)\begin{eqnarray*}
I^{\prime} & \le & \frac{N\cdot2^{N}\cdot s^{2}}{N!}\cdot\log^{N}\left(t^{N}/s\right)\\
 & \le & \frac{s^{2}\cdot e^{2N}}{N^{N}}\cdot\log^{N}\left(t^{N}/s\right)\\
 & = & \exp\left(2\log s+N\log\log t+2N-N\log\log s\right)\\
 & \le & \exp\left(2\log s+N\log\log t\right).\end{eqnarray*}
Recalling the definitions of $s$ and $t$, we finally get

\[
\log I^{\prime}\le8N\cdot\log r+C_{1}r^{2}\delta^{-2}+C_{2}r^{2}\log r\le C\left(\log r+\delta^{-2}\right)r^{2}.\]

\end{proof}
We now continue to estimate probabilities of the events $\A$ and
$\A^{\prime}$ introduced in \eqref{eq:def_set_a} and \eqref{eq:def_set_a_prime}.
\begin{lem}
\label{lem:prob_estimates}We have the following estimates:\textup{\[
\pr{\A^{\prime}\backslash\A}\le\exp\left(-\exp\left(r^{2}\right)\right)\]
and\[
\pr{\A}\le\exp\left(-S(\kappa r)+C\left(\log r+\delta^{-2}\right)r^{2}\right).\]
}\end{lem}
\begin{proof}
The first estimate in a trivial consequence of Lemma \ref{lem:log_M_upp_bnd}.
For the second we need to estimate the integral\[
I^{\prime}=\intop_{\A}\frac{1}{\pi^{N}\cdot|\Sigma|}\cdot\exp\left(-\zeta^{*}\Sigma^{-1}\zeta\right)\, d\zeta.\]
We can use the crude estimate\[
I^{\prime}\le\frac{1}{|\Sigma|}\cdot\intop_{\A}\frac{1}{\pi^{N}}\, d\zeta,\]
using the previous two lemmas, we get the result.
\end{proof}

\subsection{Lower bound for $p_{H}$}

We collect all the previous results into this
\begin{lem}
For $r$ large enough\begin{equation}
p_{H}(r)\ge S(\kappa r)-C\left(\log r+\delta^{-2}\right)r^{2},\label{eq:p_H_low_bnd_w_delta}\end{equation}
\end{lem}
\begin{proof}
Suppose that $f(z)$ has no zeros inside $\disk r$, then\[
\int_{\cir r}\log|f(z)|\, dm=\log|f(0)|.\]
We can use the fact that $\log|f(0)|$ cannot be too large, in fact\[
\pr{\log|f(0)|\ge4\log r}=\pr{|\phi_{0}|\ge r^{4}}\le\exp\left(-r^{8}\right).\]
Now using this result combined with Lemma \ref{lem:log_int_deviations}
and Lemma \ref{lem:prob_estimates}, we can bound the probability
of this event by\[
\exp\left(-r^{8}\right)+2\exp\left(-S(\kappa r)\right)+\exp\left(-S(\kappa r)+C\left(\log r+\delta^{-2}\right)r^{2}\right)\]
that is \eqref{eq:p_H_low_bnd_w_delta}.
\end{proof}
To finish the proof of the lower bound for $p_{H}(r)$, we recall
that $\kappa=1-\delta^{1/2}$ and select $\delta=r^{-\alpha}$ with
$0<\alpha<2$, then $\kappa=1-r^{-\alpha/2}$. We have to minimize
the asymptotics of the expression\[
-S(\kappa r)+C\left(\log r+\delta^{-2}\right)r^{2}.\]
Simple calculations lead to the equality $2+2\alpha=4-\alpha/2$,
or $\delta=r^{-4/5}$. We finally get (for $r$ large enough)\begin{eqnarray*}
p_{H}(r) & \ge & \frac{e^{2}}{4}\cdot\left(r\cdot\left(1-r^{-2/5}\right)\right)^{4}-Cr^{18/5}\\
 & \ge & \frac{e^{2}}{4}\cdot r^{4}-Cr^{18/5}=S(r)-Cr^{18/5}.\end{eqnarray*}

\section{Proof of The Second Theorem\label{sec:Thm_2-Proof}}

Suppose that the theorem is false, that is, there is a sequence of
entire functions\[
f_{k}(z)=\sum_{n=0}^{\infty}\phi_{n,k}\cdot\frac{z^{n}}{\sqrt{n!}},\qquad\phi_{n,k}\in K,\]
and a sequence $r_{k}\to\infty$ so that $f_{k}$ does not vanish
in $\disk{r_{k}}$. Since $K$ is a compact set, we can find a subsequence,
also denoted by $\left\{ f_{k}\right\} $, such that $\phi_{n,k}\to\phi_{n}$
for each $n\in\Bbn$. It easy to see that the sequence $\left\{ f_{k}\right\} $
converges locally uniformly to a limiting function $f$. Since $0\notin K$,
the limiting function $f$ is not identically zero. Now, using Hurwitz
theorem (see \cite[pg. 178]{Ahl}), $f$ does not vanish in any disk
$\disk{r_{k}}$; i.e. is does not vanish in the whole complex plane.

By known formulas for the order and type of entire functions by its
Taylor coefficients (see, for instance \cite[pg. 6]{Lev}) $f$ has
order $2$ and type $\frac{1}{2}$. Since it does not vanish on $\bbc$,
by Hadamard theorem, $f(z)=\exp\left(\alpha z^{2}+\beta z+\gamma\right),$
with complex constants $\alpha,\beta,\gamma$; $|\alpha|=\frac{1}{2}$.

We want to prove that we cannot get a function $f$ of this form,
using coefficients from the set $K$. We will use the asymptotics
of the coefficients of $f$ to prove this. Denoting the Taylor coefficients
of $f(z)$ by $g_{n}$, it is sufficient to show that the product\[
|g_{n}|\cdot\sqrt{n!}\]
is not bounded between any two positive constants. 

We first study the asymptotics of function of the form \eqref{eq:f_definition}.
Using Stirling's approximation we get\begin{multline}
\sqrt{n!}=\left(\sqrt{2\pi n}\left(\frac{n}{e}\right)^{n}\left(1+\O\left(\frac{1}{n}\right)\right)\right)^{1/2}\\
=\left(2\pi n\right)^{1/4}\left(\frac{n}{e}\right)^{n/2}\left(1+\O\left(\frac{1}{n}\right)\right).\label{eq:sqrt_n_fact_asymp}\end{multline}
The asymptotics of $f$ are not as simple. Using rotation and scaling,
we can assume $\alpha=\frac{1}{2}$ and $\gamma=1$, moreover it is
easy to see that $\beta$ should not be zero. Therefore the problem
is reduced to the study of the asymptotics of\begin{equation}
\exp\left(\frac{1}{2}\cdot z^{2}+\beta\cdot z\right)=\sum_{n=0}^{\infty}g_{n}\left(\beta\right)\cdot z^{n},\label{eq:limit_func_form}\end{equation}
with $\beta\ne0$, with $\beta$ possibly a complex number. A standard
application of the saddle point method shows that\begin{equation}
g_{n-1}(\beta)=C_{\beta}\cdot\left(\frac{e}{n}\right)^{\frac{n}{2}}\cdot\left(e^{\beta\sqrt{n}}+\left(-1\right)^{n}e^{-\beta\sqrt{n}}\right)\cdot\left(1+\O\left(\frac{1}{\sqrt{n}}\right)\right),\label{eq:limit_func_coeff_asymp}\end{equation}
where $C_{\beta}$ is some constant. We can see that this is not the
same rate of decay as in \eqref{eq:sqrt_n_fact_asymp} for $n\to\infty$.
We arrive at the contradiction which finishes the proof of Theorem
\ref{thm:Compact_Coeffs_Thm}.
\begin{rem*}
For the reader's convenience we prove the asymptotic estimate \eqref{eq:limit_func_coeff_asymp}
in the appendix. We also note that the function \eqref{eq:limit_func_form}
is the generating function of the Hermite polynomials.
\end{rem*}

\section{Open Problems and Further Directions}

Since it is known, that the expected number of zeros of the random
function $f(z)$ of the form \eqref{eq:f_definition}, in $\disk r$
is $r^{2}$ it is interesting to get logarithmic asymptotics for the
probability of {}``large deviations''\[
\pr{\left|n(r)-r^{2}\right|\ge\delta\cdot\sigma(r)^{\alpha}}\]
in terms of both $\alpha$ and $\delta$ (here $\sigma(r)$ denotes
the standard deviation, and $\sigma(r)=\O\left(r^{1/2}\right)$).
Some partial results in that direction can be found in Sodin and Tsirelson
\cite{ST3}, Krishnapur \cite{K} and Nazarov, Sodin and Volberg \cite{NSV}.
The most interesting cases here are $\alpha=2$ and $\alpha=4$. Note
that the aforementioned results are based on the invariance of the
random zero set and it would be also interesting to estimate similar
probabilities for more general Gaussian entire functions.

Another possible direction would be to study the hole probability
for more complicated domains. It would be interesting if there is
a function, defined for closed (simply connected) domains $U$, denote
it by $g(U)$, such that\[
\log\pr{f(z)\ne0\mbox{ in }rU}=-g(U)\cdot r^{4}+o(r^{4}).\]
Of course this function should be invariant with respect to the plane
isometries (hence {}``geometrical'').

\section{Appendix: Proof of \eqref{eq:limit_func_coeff_asymp}}

Here we compute the asymptotics of the Taylor coefficients of $\exp\left(\frac{1}{2}\cdot z^{2}+\alpha\cdot z\right)$
using the saddle point method. We write\[
\exp\left(\frac{1}{2}\cdot z^{2}+\alpha\cdot z\right)=\sum_{n=0}^{\infty}g_{n}\left(\alpha\right)\cdot z^{n}.\]
We begin with Cauchy's integral formula\begin{equation}
g_{n-1}(\alpha)=\frac{1}{2\pi i}\intop_{\disk r}\frac{\exp\left(\frac{1}{2}\cdot z^{2}+\alpha\cdot z\right)}{z^{n}}\, dz,\label{eq:g_(n-1)_def}\end{equation}
for some $r>0$ (since the function that we study is entire). We will
use the notation (notice the use of $n$ instead of $n-1$)\[
F_{n}(z)=F_{n}(\alpha;z)=\frac{1}{2}\cdot z^{2}+\alpha\cdot z-n\log z.\]
We use a standard saddle point method to study the asymptotics of
this integral. We have\[
\frac{dF_{n}}{dz}=z+\alpha-\frac{n}{z},\]
for large values of $n$ we see that the solutions to the equation
$\frac{dF_{n}}{dz}=0$ are approximately $z_{1,2}=\pm\sqrt{n}$ (notice
that they are the approximate minima of the function $F_{n}$). Therefore
we select the following rectangular contour ($T\gg1$)\begin{eqnarray*}
\Gamma & = & \left[\sqrt{n}-i\cdot T,\sqrt{n}+i\cdot T\right]\bigcup\left[-\sqrt{n}-i\cdot T,-\sqrt{n}+i\cdot T\right]\\
 &  & \bigcup\left[-\sqrt{n}-i\cdot T,\sqrt{n}-i\cdot T\right]\bigcup\left[-\sqrt{n}+i\cdot T,\sqrt{n}+i\cdot T\right]\\
 & = & \Gamma_{1}\cup\Gamma_{2}\cup\Gamma_{3}\cup\Gamma_{4}.\end{eqnarray*}
We have the following estimate for $F_{n}$ ($z=x+iy$)\begin{equation}
\left|\exp\left(F_{n}\left(z\right)\right)\right|\le\exp\left(\frac{x^{2}-y^{2}}{2}+|\alpha|\cdot|z|\right)\cdot|z|^{-n}.\label{eq:F_n_estimate}\end{equation}
We see that for $T\gg T_{0}(\alpha,n)$, $n$ fixed, we have\[
\left|\exp\left(F_{n}\left(z\right)\right)\right|\le\exp\left(-T^{2}/4\right)\]
and so\[
\left|\frac{1}{2\pi i}\intop_{\Gamma_{3}\cup\Gamma_{4}}\exp\left(F_{n}(z)\right)\, dz\right|\le C\sqrt{n}\cdot\exp\left(-T^{2}/4\right)\rightarrow0,\]
as $T\rightarrow\infty$. Therefore what is left is to estimate the
integrals\begin{eqnarray*}
I_{1} & = & \frac{1}{2\pi i}\cdot\intop_{\sqrt{n}+i\Bbr}\exp\left(F_{n}(z)\right)\, dz,\\
I_{2} & = & \frac{1}{2\pi i}\cdot\intop_{-\sqrt{n}-i\Bbr}\exp\left(F_{n}(z)\right)\, dz.\end{eqnarray*}
We start by using the parametrization $z=\sqrt{n}\cdot\left(1+i\cdot t\right)$,
with $t\in\left(-\infty,\infty\right)$. Then we have\begin{eqnarray*}
I_{1} & = & \frac{n^{1/2}}{2\pi}\cdot\intop_{-\infty}^{\infty}\exp\left(\frac{n}{2}\cdot\left(1+i\cdot t\right)^{2}+\alpha\sqrt{n}\cdot\left(1+i\cdot t\right)\right)\cdot n^{-n/2}\cdot\left(1+i\cdot t\right)^{-n}\, dt\\
 & = & \frac{1}{2\pi}\cdot\frac{\exp\left(\frac{n}{2}+\alpha\sqrt{n}\right)}{n^{\left(n-1\right)/2}}\cdot\intop_{-\infty}^{\infty}\exp\left(-\frac{n}{2}\cdot t^{2}+\left(n+\alpha\sqrt{n}\right)\cdot ti\right)\cdot\left(1+i\cdot t\right)^{-n}\, dt.\end{eqnarray*}
We denote the new integral by $I_{1}^{\prime}.$ Using \eqref{eq:F_n_estimate}
we want to show that only the {}``small'' values of $t$ contribute
to the asymptotics. 

Set $0<b<\frac{1}{2}$, we want $b$ to be sufficiently close to $\frac{1}{2}$
and to have\[
-\frac{n}{2}\cdot t^{2}+\sqrt{n}\left|\alpha\right|\cdot|t|\le-b\cdot nt^{2},\]
or after rearrangement\[
|t|\ge\frac{|\alpha|}{\sqrt{n}}\cdot\frac{1}{\frac{1}{2}-b}.\]
We assume $n\ge2^{24}$ and select $b=\frac{1}{2}-\frac{1}{n^{1/12}}$
(the reason for selecting $\frac{1}{12}$ will be clear from the calculations
bellow). For $|t|\ge\frac{|\alpha|}{n^{5/12}}$, we now have\begin{eqnarray*}
\left|\exp\left(-\frac{n}{2}\cdot t^{2}+\left(n+\alpha\sqrt{n}\right)\cdot ti\right)\right| & \le & \exp\left(-\frac{n}{2}\cdot t^{2}+\sqrt{n}\left|\alpha\right|\cdot|t|\right)\\
 & \le & \exp\left(-\left(\frac{1}{2}-\frac{1}{n^{1/12}}\right)\cdot nt^{2}\right)\\
 & \le & \exp\left(-\frac{1}{4}\cdot nt^{2}\right).\end{eqnarray*}
We also note that for $n$ large enough\[
\left|1+i\cdot\frac{|\alpha|}{n^{5/12}}\right|^{n}\ge\left(1+\frac{|\alpha|^{2}}{n^{5/6}}\right)^{n/2}\ge\exp\left(\frac{1}{3}\cdot|\alpha|^{2}n^{1/6}\right).\]
Therefore we have the following estimate for {}``large'' values
of $t$\begin{eqnarray*}
\left|\intop_{|t|\ge\frac{|\alpha|}{n^{5/12}}}\exp\left(-\frac{n}{2}\cdot t^{2}+\left(n+\alpha\sqrt{n}\right)\cdot ti\right)\cdot\left(1+i\cdot t\right)^{-n}\, dt\right| & \le\\
\intop_{|t|\ge\frac{|\alpha|}{n^{5/12}}}\exp\left(-\frac{n}{4}\cdot t^{2}\right)\, dt\cdot\left|1+i\cdot\frac{|\alpha|}{n^{5/12}}\right|^{-n} & \le\\
\intop_{-\infty}^{\infty}\exp\left(-\frac{n}{4}\cdot t^{2}\right)\, dt\cdot\exp\left(-\frac{1}{3}\cdot|\alpha|^{2}n^{1/6}\right) & \le\\
\exp\left(-\frac{1}{3}\cdot|\alpha|^{2}n^{1/6}\right).\end{eqnarray*}
For the main part we use the Taylor expansion of $\log\left(1+x\right)$
for small $x$, \[
\log\left(1+i\cdot t\right)=i\cdot t+\frac{t^{2}}{2}+\sum_{m=3}^{\infty}\left(-1\right)^{m+1}\cdot\frac{\left(it\right)^{m}}{m}.\]
For $|t|<\frac{|\alpha|}{n^{5/12}}$ and $n^{5/12}\ge2|\alpha|$,
we have the following bound for the tail\[
\left|\sum_{m=3}^{\infty}\left(-1\right)^{m+1}\cdot\frac{\left(it\right)^{m}}{m}\right|\le|t|^{3}\cdot\left|\sum_{m=0}^{\infty}\left(\frac{1}{2}\right)^{m}\right|=2|t|^{3},\]
so we get\[
\left|\log\left(1+i\cdot t\right)-\left(i\cdot t+\frac{t^{2}}{2}\right)\right|\le2|t|^{3},\]
or\[
\left|n\cdot\log\left(1+i\cdot t\right)-n\left(i\cdot t+\frac{t^{2}}{2}\right)\right|\le2n|t|^{3}.\]
Using $\exp\left(n|t|^{3}\right)=1+n|t|^{3}+\O\left(\frac{1}{\sqrt{n}}\right)$,
we now have\[
I_{1}^{\prime}=I_{1}^{\dprime}+I_{1}^{\left(1\right)}+\O\left(\frac{1}{\sqrt{n}}\cdot I_{1}^{\left(2\right)}\right).\]
where \begin{eqnarray*}
I_{1}^{\dprime} & = & \intop_{|t|<\frac{|\alpha|}{n^{5/12}}}\exp\left(-nt^{2}+\alpha\sqrt{n}\cdot ti\right)\, dt,\\
I_{1}^{\left(1\right)} & = & \intop_{|t|<\frac{|\alpha|}{n^{5/12}}}\exp\left(-nt^{2}+\alpha\sqrt{n}\cdot ti\right)\cdot2n|t|^{3}\, dt,\\
I_{1}^{\left(2\right)} & = & \intop_{|t|<\frac{|\alpha|}{n^{5/12}}}\exp\left(-nt^{2}\right)\, dt.\end{eqnarray*}
After simple approximations of $I_{1}^{\left(1\right)}$ and $I_{1}^{\left(2\right)}$
we get\[
I_{1}^{\prime}=I_{1}^{\dprime}+\O\left(\frac{1}{n}\right).\]
Using similar estimates we conclude\[
\left|I_{1}^{\dprime}-\intop_{-\infty}^{\infty}\exp\left(-nt^{2}+\alpha\sqrt{n}\cdot ti\right)\, dt\right|\le\intop_{|t|\ge\frac{|\alpha|}{n^{5/12}}}\exp\left(-\frac{n}{2}\cdot t^{2}\right)\, dt,\]
the error term is bounded in the following way\begin{eqnarray*}
\intop_{|t|\ge\frac{|\alpha|}{n^{5/12}}}\exp\left(-\frac{n}{2}\cdot t^{2}\right)\, dt & \le & \intop_{|t|\ge\frac{|\alpha|}{n^{5/12}}}\exp\left(-\frac{n}{2}\cdot t^{2}\right)\cdot|t|\cdot\frac{n^{5/12}}{|\alpha|}\, dt\\
 & = & \frac{2\exp\left(-\frac{|\alpha|^{2}}{2}\cdot n^{1/6}\right)}{|\alpha|n^{7/12}}.\end{eqnarray*}
Overall we have (for $n$ large enough)\[
I_{1}^{\prime}=\intop_{-\infty}^{\infty}\exp\left(-nt^{2}+\alpha\sqrt{n}\cdot ti\right)\, dt+\O\left(\frac{1}{n}\right)=\sqrt{\frac{\pi}{n}}\cdot\exp\left(-\alpha^{2}/4\right)+\O\left(\frac{1}{n}\right),\]
and so\[
I_{1}=\sqrt{\frac{1}{4\pi}}\cdot\frac{\exp\left(\frac{n}{2}+\alpha\sqrt{n}-\alpha^{2}/4\right)}{n^{n/2}}\cdot\left(1+\O\left(\frac{1}{\sqrt{n}}\right)\right).\]
Repeating the same calculation for $I_{2}$ we also have\[
I_{2}=\sqrt{\frac{1}{4\pi}}\cdot\left(-1\right)^{n}\cdot\frac{\exp\left(\frac{n}{2}-\alpha\sqrt{n}-\alpha^{2}/4\right)}{n^{n/2}}\cdot\left(1+\O\left(\frac{1}{\sqrt{n}}\right)\right).\]
Recalling \eqref{eq:g_(n-1)_def}\[
g_{n-1}(\alpha)=C_{\alpha}\cdot\left(\frac{e}{n}\right)^{\frac{n}{2}}\cdot\left(e^{\alpha\sqrt{n}}+\left(-1\right)^{n}e^{-\alpha\sqrt{n}}\right)\cdot\left(1+\O\left(\frac{1}{\sqrt{n}}\right)\right).\]
Comparing this with \eqref{eq:sqrt_n_fact_asymp} we get the required
contradiction.

\end{document}